\documentclass[12pt]{article}

\usepackage{latexsym}
\usepackage{amsfonts}
\usepackage{amsmath}
\usepackage{amssymb}

\newcommand\A{{\mathbb A}}
\newcommand\Z{{\mathbb Z}}
\newcommand\Q{{\mathbb Q}}
\newcommand\C{{\mathbb C}}
\newcommand\xh{{\hat{x}}}
\newcommand\uAut{{\underline{\operatorname{Aut}}}}
\newcommand\spec{{\mathsf{Spec}\,}}
\newcommand\Gm{{{\mathbb G}_m}}
\newcommand\D{{\mathcal D}}

\renewcommand\k{{\bf k}}
\newcommand\OO{{\mathcal O}}

\newtheorem{lmm}{Lemma}
\newtheorem{thm}{Theorem}
\newtheorem{prop}{Proposition}
\newtheorem{conj}{Conjecture}
\newtheorem{rmk}{Remark}

\def\ch{\operatorname{char}}
\def\Ext{\operatorname{Ext}}
\def\Center{\operatorname{Center}}
\def\Der{\operatorname{Der}}
\def\Ker{\operatorname{Ker}}
\def\Lie{\operatorname{Lie}}
\def\Id{\operatorname{Id}}
\def\Br{\operatorname{Br}}
\def\Sp{\operatorname{Sp}}
\def\Fr{\operatorname{Fr}}
\def\Mat{\operatorname{Mat}}
\def\ad{\operatorname{ad}}
\def\Aut{\operatorname{Aut}}

\title{Automorphisms of the Weyl algebra}

\begin{document}

\author{Alexei Belov-Kanel  and Maxim Kontsevich}
\maketitle

\vspace{5mm}

\vspace{5mm}
\noindent\textbf{Abstract.}
 We discuss a conjecture which says that the automorphism group
 of the Weyl algebra in characteristic zero is canonically isomorphic
 to the automorphism group of the corresponding Poisson algebra
 of classical polynomial symbols. Several arguments in favor
 of this conjecture are presented, all based on the consideration of the
 reduction of the Weyl algebra to positive characteristic.

\vspace{5mm}

\section{Introduction}\label{sec:conj}

This paper is devoted to the following surprising conjecture.

\begin{conj}\label{conj:iso}
The automorphism group of the Weyl algebra of index $n$ over $\C$ is
isomorphic to the group of the polynomial symplectomorphisms of
a $2n$-dimensional affine space
 $$\Aut(A_{n,\C})\simeq
 \Aut(P_{n,\C})\,\,.$$
 \end{conj}

Here for an integer $n\ge 1$, denote by $A_{n,\C}$ the Weyl algebra of
index $n$ over $\C$
$$
 \C\langle \xh_1,\dots ,\xh_{2n}\rangle/\left(\mbox{ relations
 }[\xh_i,\xh_j]=\omega_{ij}\,,\,\, \,1\le i,j\le  2n\,\right)
$$
where $(\omega_{ij})_{1\le i,j\le 2n}$ is the standard skew-symmetric
matrix:
$$\omega_{ij}=\delta_{i,n+j}-\delta_{n+i,j}\,\,,$$
and by $P_{n,\C}$ the Poisson algebra over $\C$ which is the usual
polynomial algebra $\C[x_1,\dots,x_{2n}]\simeq \OO(\A^{2n}_\C)$ endowed
with the Poisson bracket:
$$ \{x_i,x_j\}=\omega_{ij}\,,\,\,1\le i,j\le
2n\,\,.$$

The algebra $A_{n,\C}$ is isomorphic to the algebra $\D(\A^n_\C)$ of
polynomial differential operators in $n$ variables $x_1,\dots,x_n$:

$$
\xh_i\mapsto x_i,\,\,\,\xh_{n+i}\mapsto \partial_i:=\partial/\partial
x_i,\,\,1\le i\le n\,\,.
$$

The conjecture becomes even more surprising if one takes into account
the fact that the {\it Lie algebras} of derivations of $A_{n,\C}$ and
 $P_{n,\C}$ are not isomorphic to each other (see Section
\ref{sec:negative}). Conjecture \ref{conj:iso} is closely related to a
question raised by one of us few years ago   (see section 4.1,
Question 4 in \cite{Kontsevich1}), the motivation at that time came
from the theory of deformation quantization for algebraic varieties.

One of  main results of the present paper is Theorem \ref{thm:tame}
which says that the subgroups of $\Aut(A_{n,\C})$ and $\Aut(P_{n,\C})$
consisting of the so-called {\it tame} automorphisms, are naturally
isomorphic to each other. Another result is Theorem \ref{thm:frob}
from Section \ref{sec:cand}, which allows to propose a (hypothetical)
specific  candidate for the isomorphism between two automorphism groups
as above. The key idea is to use reduction to finite characteristic and
the fact that the Weyl algebra in finite characteristic has a huge
center isomorphic to the polynomial algebra. Another application of
this idea is a proof of a stable equivalence between the Jacobian
and Dixmier conjectures, see \cite{JD}.  This paper is an
extended version of  the talk given by one of us (see
\cite{Kontsevich2}) on Arbeitstagung 2005 (Bonn).

We finish the introduction with

     \subsection{First positive evidence: case $n=1$}

The structure of the group $\Aut(P_{1,\C})$ is known since the work of H.~W.~E.~Jung
(see \cite{Jung}).  This group contains the group $G_1=SL(2,\C)\ltimes
\C^2$ of special affine transformations, and the solvable group $G_2$ of
polynomial transformations of the form
$$
(x_1,x_2)\mapsto (\lambda x_1+ F(x_2),\lambda^{-1} x_2),\,\,\lambda\in
\C^\times,\,\, F\in \C[x]\,\,.
$$

The group $\Aut(P_{1,\C})$ is equal to the  amalgamated product of $G_1$
and $G_2$ over their intersection. J.~Dixmier in \cite{Dixmier} and
later L.~Makar-Limanov in \cite{Makar-Limanov} proved that  if one
replaces the commuting variables $x_1,x_2$ by noncommuting variables
$\xh_1,\xh_2$ in the formulas above, one obtains the description of the group
$\Aut(A_{1,\C})$.  Hence, in the case $n=1$ the two automorphism groups are
 isomorphic.

\bigskip

{\bf Acknowledgments}: We are grateful to Ofer Gabber, Leonid
Makar-Limanov and Eliyahu Rips for fruitful discussions and remarks.

  \section{Automorphism groups as ind-schemes}

For an arbitrary commutative ring $R$ one can define the Weyl algebra $A_{n,R}$ over $R$,
by just replacing $\C$ by $R$ in the definition. We denote the algebra $A_{n,\Z}$  simply
by $A_n$, hence $A_{n,R}=A_n\otimes R$.
The algebra $A_{n,R}$ considered as an $R$-module is free with basis
$$\xh^\alpha:=\xh_1^{\alpha_1}\dots \xh_{2n}^{\alpha_{2n}},\,\,\alpha=
(\alpha_1,\dots,\alpha_{2n})\in\Z_{\ge 0}^{2n}\,\,.$$

We define an increasing filtration (the Bernstein filtration) on the algebra
$A_{n,R}$ by
$$A_{n,R}^{\le N}:=\left\{\sum_{\alpha} c_\alpha \xh^\alpha\,|
 \,c_\alpha\in R,\,\,c_\alpha=0\mbox{ for
 }|\alpha|:=\alpha_1+\dots+\alpha_{2n}>N\right\}\,\,.$$
This filtration induces a filtration on the automorphism group:
$$\Aut^{\le N} A_{n,R}:=\{f\in \Aut(A_{n,R})| \,f(\xh_i),f^{-1}(\xh_i)\in A^{\le N}_{n,R}\,\,\,\forall i=1,\dots,2n\}\,\,.$$

The following functor on commutative rings:
  $$R\mapsto \Aut^{\le N} (A_{n,R})\,\,,$$
is representable by an affine scheme of finite type over $\Z$. We
denote this scheme by $$\uAut^{\le N}(A_n)\,\,.$$ The ring of functions
$\OO(\uAut^{\le N}(A_n))$ is generated by the variables
$$(c_{i,\alpha},c'_{i,\alpha})_{1\le i\le 2n,|\alpha|\le N}$$
which are the coefficients of the elements $f(\xh_i),f^{-1}(\xh_i)$ in the
standard basis $(\xh^\alpha)$ of the Weyl algebra.

The obvious inclusions $\uAut^{\le N}(A_n)\hookrightarrow\uAut^{\le
(N+1)}(A_n)$ are closed embeddings, the inductive limit over $N$ of
schemes $\uAut^{\le N}(A_n)$ is an ind-affine scheme over $\Z$.  We
denote it by $\uAut(A_n)$, it is a group-like object in the category
of ind-affine schemes.

Similarly, one can define all the above notions for the Poisson
algebra $P_n$, in particular we  have an affine scheme
$\uAut^{\le N}(P_n)$ of finite type, and a group ind-affine
scheme $\uAut(P_n)$.

Later it will be convenient to use the notation $$\uAut^{\le
N}(A_{n,R}):=\uAut^{\le N}(A_n)\times_{\spec \Z} \spec R$$ (here $R$ is
an arbitrary commutative ring), for a scheme over $\spec R$ obtained by
the extension of scalars, similarly we will have schemes $\uAut^{\le
 N}(P_{n,R})$ for the case of Poisson algebras.

There is also another sequence of closed embeddings (stabilization)
$$\uAut^{\le N}(A_n)\hookrightarrow\uAut^{\le N}(A_{n+1}),\,\,
\uAut^{\le N}(P_n)\hookrightarrow\uAut^{\le N}(P_{n+1})$$ corresponding
to the addition of two new generators and extending the automorphisms
by the trivial automorphism on the additional generators.

Conjecture \ref{conj:iso} says that groups of points $\uAut(A_n)(\C)$
and $\uAut(P_n)(\C)$ are isomorphic.  We expect that the isomorphism
should preserve the filtration by degrees, compatible with
stabilization embeddings, and should be a constructible map for any
given term of filtration, defined over $\Q$:  \begin{conj}
\label{conj:const} There exists a family $\phi_{n,N}$ of constructible
one-to-one maps $$\phi_{n,N}:\uAut^{\le N}(A_{n,\Q}) \to \uAut^{\le
N}(P_{n,\Q})$$ compatible with the  inclusions increasing  indices $N$
and $n$, and with the group structure.  \end{conj}

  Obviously, this conjecture implies Conjecture \ref{conj:iso}, and
  moreover it implies that the isomorphism exists if one replaces $\C$
by an arbitrary field of characteristic zero.

 	   \section{Negative evidence}\label{sec:negative}

For any group object $G$ in the category of ind-affine schemes over
$\Q$, one can associate its Lie algebra $\Lie (G)$, by considering
points of $G$ with coefficients in the algebra of dual numbers
$\Q[t]/(t^2)$.  The Lie algebras
 
$$\Lie(\uAut(A_{n,\Q})),\,\,\,\Lie(\uAut(P_{n,\Q}))$$
are by definition the algebras of derivations of $A_{n,\Q}$ and
$P_{n,\Q}$, respectively.  The following fact is well known:

\begin{lmm} All derivations of the Weyl algebra are inner:
$$\Der(A_{n,\Q})\simeq A_{n,\Q}/\Q\cdot 1_{A_{n,\Q}},\,\,\,f\in A_{n,\Q}\mapsto
[f,\cdot\,]\in \Der(A_{n,\Q})\,\,.$$
\end{lmm}

{\bf Proof}: For any algebra $A$ the space of its outer derivations 
coincides with the first Hochschild cohomology,
 $$ H^1(A,A)=\Ext^1_{A-mod-A}(A,A)\,\,. $$

We claim that for $A=A_{n,\Q}$ the whole Hochschild cohomology is
a 1-dimensional space in degree zero. Namely, it is easy to see that
there exists an isomorphism $A\otimes A^{op}\simeq \D(\A^{2n}_\Q)$ such
that the diagonal bimodule $A$ is isomorphic to $\OO(\A^{2n}_\Q)$. The
result follows from a standard property of $\D$-modules:
 $$\Ext^*_{\D_X-mod}(\OO_X,\OO_X)=H^*_{de\,\,Rham}(X)\,\,$$ which holds
for any smooth variety $X$.  $\Box$

Derivations of $P_{n,\Q}$ are polynomial Hamiltonian vector fields:

$$\Der(P_{n,\Q})\simeq P_{n,\Q}/\Q\cdot 1_{P_{n,\Q}},\,\,\,f\in
P_{n,\Q}\mapsto \{f,\cdot\, \}\in \Der(P_{n,\Q})\,\,.$$

We see that both Lie algebras of derivations are of the ``same'' size,
each of them has a basis labeled by the set $\Z_{\ge 0}^{2n}\setminus
 \{(0,\dots,0)\}$. Nevertheless, these two Lie algebras are not
 isomorphic. Namely, $\Der(P_{n,\Q})$ has many nontrivial Lie
subalgebras of finite codimension (e.g.  Hamiltonian vector fields
vanishing at a given point in $\Q^{2n}$), whereas the algebra
$\Der(A_{n,\Q})$ has no such subalgebras.

The conclusion is that the hypothetical constructible isomorphism
$\phi_{n,N}$ cannot be a scheme map.

 \section{Positive evidence: tame automorphisms}\label{sec:tame}

  Obviously, the symplectic group $\Sp(2n,\C)$ acts by automorphisms of
$A_{n,\C}$ and  $P_{n,\C}$, by symplectic   linear transformations of the 
variables  $\xh_1,\dots,\xh_{2n}$ and  $x_1,\dots,x_{2n}$,
respectively.  Also, for any polynomial $F\in \C[x_1,\dots,x_n]$ we
define   ``transvections''

$$T_F^A\in \Aut(A_{n,\C}),\,\, T_F^P\in \Aut(P_{n,\C})$$
by the formulas
$$T_F^P(x_i)=x_i,\,\, T_F^P(x_{n+i})=x_{n+i}+\partial_i
F(x_1,\dots,x_n),\,\,1\le i\le n\,\,,$$
$$T_F^A(\xh_i)=\xh_i,\,\,T_F^A(\xh_{n+i})=\xh_{n+i}+\partial_i F
(\xh_1,\dots,\xh_n),\,\,1\le i\le n\,\,.$$
The last formula makes sense, as we substitute the commuting variables
$\xh_1,\dots,\xh_n$ in place of $x_1,\dots,x_n$ in the polynomial
$\partial_i F:=\partial F/\partial x_i$. A straightforward check shows
that these maps are well defined, i.e. $T_F^P$ preserves the Poisson
bracket and $T_F^A$ preserves the commutation relations between the generators.
The automorphism $T_F^A$ is in a sense the conjugation by a nonalgebraic
element $\exp(F(\xh_1,\dots,\xh_n))$.

The correspondence $F\mapsto T_F^A$ (resp.  $F\mapsto T_F^P$) gives a group
homomorphism $\C[x_1,\dots,x_n]/\C\cdot 1\to \Aut(A_{n,\C})$ (resp. to
$\Aut(P_{n,\C})$).

Let us denote by $G_n$ the free product of $\Sp(2n,\C)$ with the abelian
group $\C[x_1,\dots,x_n]/\C\cdot 1$.  We obtain two homomorphisms
$\rho^A_n$ and $\rho^P_n$ from $G_n$ to $\Aut(A_{n,\C})$ and
$\Aut(P_{n,\C})$ respectively. The automorphisms which belong to the image
of $\rho^A_n$ (resp. $\rho^P_n$) are called {\it tame}.

 \begin{thm} \label{thm:tame}
In the above notation, $\Ker \rho^A_n=\Ker \rho^P_n$.
\end{thm}

As an immediate corollary, we obtain that the groups of tame
automorphisms of $A_{n,\C}$ and  $P_{n,\C}$ are canonically
isomorphic.  It is quite possible that all of the automorphisms of $A_{n,\C}$
and $P_{n,\C}$ become tame after stabilization, i.e. after adding
several dummy variables and increasing the parameter $n$.  If this is the
case, then we obtain Conjecture \ref{conj:iso} (and, in fact,
Conjecture \ref{conj:const} as well).

\begin{rmk} We expect that the canonical isomorphism in Conjecture
\ref{conj:const} coincides with the above isomorphism on subgroups of
tame elements.
\end{rmk}

The proof of  Theorem \ref{thm:tame} will be given in Section
\ref{sec:tame2}.

\begin{rmk} 

I.R.Shafarevich in \cite{Shafarevich1} introduced a notion of a Lie algebra
 for an infinite-dimensional algebraic group (in fact, his notion of a group is a bit obscure as he 
 does not use the language of  ind-schemes).
  It is known that the Lie algebra associated with the group of all of the
automorphisms of a polynomial ring coincides with the Lie algebra of the group
of {\rm tame} automorphisms.  The same is true for the automorphisms of the Weyl
algebra and polynomial symplectomorphisms.
However, I.P.Shestakov and U.Umirbaev \cite{ShUm} have shown that
the Nagata automorphism is wild. Hence, we have an infinite-dimensional
effect: A proper subgroup  can have the same Lie
algebra as the whole group. Our Conjecture indicates that further pathologies
 are possible, the same group has two different Lie algebras when interpreted
  as an ind-scheme in two different ways. Presumably, it means that at least
  one of our automorphism groups is singular everywhere.
\end{rmk}

    \section{The Weyl algebra in finite characteristic}

    Here we introduce the main tool which allows us to relate  the algebras
$A_{n,R}$ and $P_{n,R}$ (and their automorphisms).

     \subsection{The Weyl algebra as an Azumaya algebra}

     Let $R$ be a commutative ring in characteristic $p>0$, i.e.
$p\cdot 1_R=0\in R$.  It is well known that in this case the Weyl
algebra $A_{n,R}$ has a big center, and moreover, it is an Azumaya
algebra of rank $p^n$ (see \cite{Revoy}).

\begin{prop}
For any commutative ring $R\supset \Z/ p\Z$ the center
$C_{n,R}$ of $A_{n,R}$ is isomorphic as an $R$-algebra to the
polynomial algebra $R[y_1,\dots,y_{2n}]$, where the variable
$y_i,\,\,i=1,\dots,2n$, corresponds to $\xh_i^p$. The algebra
 $A_{n,R}$ is a free $C_{n,R}$-module of rank $p^{2n}$, and
it is an Azumaya algebra of rank $p^n$ over $C_{n,R}$.
\end{prop}

{\bf Proof}: First of all, a straightforward check shows that the elements
$(\xh_i^p)_{1\le i\le 2n}$ are central and generate over $R$ the
polynomial algebra. The algebra $A_{n,R}$ is the algebra over
$R[y_1,\dots,y_{2n}]$ with generators $\xh_1,\dots,\xh_{2n}$ and
relations

$$[\xh_i,\xh_j]=\omega_{ij},\,\,\xh_i^p=y_i,\,\,1\le i,j\le 2n\,\,.$$

After the extension of scalars from
$R[y_1,\dots,y_{2n}]$ to
$C'_{n,R}:=R[y_1^{1/p},\dots,y_{2n}^{1/p}]$ the pullback of the
Weyl algebra can be described as an algebra over $C'_{n,R}$
with generators $\xh_i':=\xh_i-y_i^{1/p}$ and relations
$$[\xh_i',\xh_j']=\omega_{ij},\,\,(\xh_i')^p=0,\,\,1\le i,j\le 2n\,\,.$$

It is well known that the algebra over $\Z/p \Z$ with two generators
$\xh_1,\xh_2$ and defining relations
$[\xh_1,\xh_2]=1,\,\xh_1^p=\xh_2^p=0$ is isomorphic to the matrix
algebra $\Mat(p\times p,\Z/p\Z)$ (consider operators $x$ and $d/dx$ in
the truncated polynomial ring $\Z/p \Z[x]/(x^p)$). Hence, after the
faithfully flat finitely generated extension from $R[y_1,\dots,y_{2n}]$
to $C'_{n,R}$, we obtain the matrix algebra $\Mat(p^n\times p^n,C'_{n,R})$.
 Then the proposition follows from standard properties of Azumaya
 algebras. $\Box$

     \subsection{Poisson bracket on the center of the Weyl algebra}

The next observation is that for any commutative ring $R$ flat over
any prime $p\in \spec \Z$ one can define a Poisson bracket on $C_{n,R/pR}$
in an intrinsic manner.  Namely, for such $R$ and any two elements
$a,b\in C_{n,R/pR}$ we define the element $\{a,b\}\in C_{n,R/pR}$ by the
formula
$$
\{a,b\}:=\frac {[\tilde{a},\tilde{b}]}{p} \pmod{p R}\,\,,
$$
where $\tilde{a},\tilde{b}\in A_{n,R}$ are arbitrary lifts of

$$a,b\in  C_{n,R/pR}\subset A_{n,R/pR}= A_{n,R}\pmod{p R}\,\,.$$

 A straightforward check shows that the operation $(a,b)\mapsto
\{a,b\}$ is well defined, takes values in $C_{n,R/pR}\subset
A_{n,R/pR}$, satisfies the Leibniz rule with respect to the product
on $C_{n,R/pR}$, and the Jacobi identity.

Morally, our construction of the bracket is analogous to the well-known
counterpart in deformation quantization.  If one has a
one-parameter family of associative algebras $A_{\hbar}$ (flat over the
algebra of functions of $\hbar$), then on the center of $A_0$ one has a
canonical Poisson bracket given by the ``same'' formula as above:
$$
\{a_0,b_0\}:=\frac{[a_\hbar,b_\hbar]}{\hbar}+O(\hbar)\,\,,
$$
where $a_\hbar,b_\hbar\in A_\hbar$ are arbitrary extensions of elements
$$a_0,b_0\in \Center(A_0)\subset A_0\,\,.$$

The prime number $p$ plays the role of Planck constant $\hbar$.

 The following lemma shows that the canonical Poisson bracket on
$C_{n,R/pR}$ coincides (up to sign) with the standard one:

  \begin{lmm} In the above notation, one has
  $$\{y_i,y_j\}=-\omega_{ij}\,\,.$$
  \end{lmm}

{\bf Proof}: It is enough to make the calculation in the case of one variable.  The
following elementary identity holds in the algebra of differential
operators with coefficients in $\Z$:
$$
[(d/dx)^p,x^p]=\sum_{i=0}^{p-1} \frac{(p!)^2}{(i!)^2 (p-i)!} x^i
(d/dx)^i\,\,.
$$
The r.h.s. is divisible by $p$, and is equal to $-p$ modulo $p^2$.
$\Box$

In the above considerations one can make a weaker assumption,
 it is enough to consider the coefficient ring $R$ flat over $\Z/p^2 \Z$.
 The corollary is that for any automorphism of the Weyl algebra in characteristic $p$
  which admits a lift$\mod {p^2}$, the induced automorphism
   of the center preserves the canonical Poisson bracket.
    The condition of the existence of the lift is necessary. For example,
    the automorphism of $A_{2,\Z/p\Z}$ given by
   $$\xh_1\mapsto \xh_1+\xh_2^p\xh_3^{p-1},\,\,\,\,\,\xh_i\mapsto \xh_i,\,i=2,3,4$$
   acts on the center by
   $$y_1\mapsto y_1+y_2^p y_3^{p-1}-y_2,\,\,\,\,\,y_i\mapsto y_i,\,i=2,3,4\,\,.$$
   The above map does not preserve the Poisson bracket.

 \section{Correspondence between automorphisms in finite characteristic}
 \subsection{Rings at infinite prime} \label{sec:rinfty}
  It will be convenient to introduce the following notation
(``reduction modulo infinite prime'') for an arbitrary commutative ring
  $R$:
  $$
  R_\infty:=\lim_{\stackrel{\longrightarrow}{f.g.\,R'\subset
  R}}\left(\,\,\prod_{\mathrm{primes}\,\,p} R'\otimes \Z/p \Z
  \,\,\,\,\big{/} \, \bigoplus_{\mathrm{primes}\,\,p} R'\otimes \Z/
  p\Z\right)\,\,.
  $$

  Here the inductive limit is taken over the filtered system consisting
  of all finitely generated subrings of $R$, the index $p$ runs over
  primes $2,3,5,\dots$.

It is easy to see that the ring $R_\infty$ is defined over $\Q$ (all
primes are invertible in $R_\infty$), and the obvious map $R\mapsto
R_\infty$ gives an inclusion $i:R\otimes \Q\hookrightarrow R_\infty$.
Also, there is a universal Frobenius endomorphism $\Fr:R_\infty\to
R_\infty$ given by

$$\Fr(a_p)_{\,\mathrm{primes}\,\,p}:=(a_p^p)_{\,\mathrm{primes}\,\,p}\,\,.$$

Finally, if $R$ has no nilpotents then $\Fr\circ i$ gives another inclusion of $R\otimes \Q$ into $R_\infty$.

\begin{rmk} Maybe a better notation would be $R\pmod{\infty}$ instead of $R_\infty$ as one can also imitate
 $p$-adic completion:

 $$
\lim_{\stackrel{\longrightarrow}{f.g.\,R'\subset R}}\,\,\lim_{\stackrel{\longleftarrow}{n\ge 1}}
 \left(\,\,\prod_{\mathrm{primes}\,\,p} R'\otimes \Z/p^n \Z \,\,\,\,\big{/} \,
  \bigoplus_{\mathrm{primes}\,\,p} R'\otimes \Z/ p^n\Z\right)\,\,.
$$
This larger ring has a canonical element ``infinite prime'' $P$ which
is the class of sequence $a_p=p\,\,\,\forall$ prime $p$. Our
``reduction modulo infinite prime'' is literally the reduction of the
larger ring modulo $P$.
\end{rmk}

	 \subsection{The homomorphism $\psi_R$}

It follows from the previous section that one has a canonical group
homomorphism
 $$\psi_R:\uAut(A_n)(R)\to \uAut(P_n)(R_\infty)\,\,.$$
Namely, if $f\in \uAut(A_n)(R)$ is an automorphism of $A_{n,R}$ then it
belongs to a certain term of filtration $\uAut^{\le N}(A_n)(R)$ and 
moreover, is defined over a finitely generated ring $R'\subset R$. For any
prime $p$, the automorphism $f$ gives an automorphism $f_p$ of $A_{n,R'/p
R'}$, and hence an automorphism $f_p^{centr}$ of the center $C_{n,R'/p
 R'}\simeq R'/p R'[y_1,\dots,y_{2n}]$.

\begin{lmm} \label{lmm:deg}
For any $f\in \Aut^{\le N}(A_n)(R)$ and any $i=1,\dots,2n$,
  the element $f_p^{centr}(y_i)\in C_{n,R'/pR'}$ is a polynomial
 in $y_1,\dots,y_{2n}$ of degree $\le N$.
 \end{lmm}

 {\bf Proof}:  One has $f_p^{centr}(y_i)=f_p(\xh_i^p)=(f_p(\xh_i))^p$.
The element $f_p(\xh_i)\in A_{n,R'/p R'}$ has degree $\le N$ by our
assumption. Hence, $f_p(\xh_i^p)$ has degree $\le pN$.  We know that
the last element is in fact a polynomial in the commuting variables
$\xh_1^p,\dots,\xh_{2n}^p$. Then it is a polynomial of degree $\le N$
in these variables, as it follows immediately from the fact that
$(\xh^\alpha)_{\alpha\in\Z_{\ge 0}^{2n}}$ is a $R'/pR'$ basis of
$A_{n,R'/p R'}$. $\Box$

Any finitely generated commutative ring  is flat over all sufficiently
large primes. Hence, we obtain polynomial symplectomorphisms
$f_p^{centr}$ for $p\gg 1$ of degree (and the degree of the inverse
automorphism) uniformly bounded by $N$ from above.  We define
$\psi_R(f)$ to be the collection $(f_p^{centr})_{p\gg 1}$ of
automorphisms of $P_{n,R'/p R'}$ where we identify the variables $y_i\in
C_{n,R'/p R'}$ with $x_i\in P_{n,R'/p R'}$, considered asymptotically
in $p$.  It is easy to see that $\psi_R(f)$ is an element of
$\uAut^{\le N}(P_n)(R_\infty)$ and it does not depend on the choice of
a finitely generated ring $R'\subset R$ over which $f$ is defined.
Hence we obtain a canonical map $\psi_R$ which is obviously a group
homomorphism.

\subsection{Untwisting the Frobenius endomorphism}\label{sec:cand}

	       \begin{thm}	       \label{thm:frob}
Let $R$ be a finitely generated ring such that $\spec R$ is smooth over
$\spec \Z$.  Then for any $f\in \Aut(A_{n,R})$ the corresponding
symplectomorphism $f_p^{centr}\in \Aut(P_{n,R/pR})$ is defined over
 $(R/pR)^p$ for sufficiently large $p$.
 \end{thm}

{\bf Proof}:  First of all, notice that the subring $(R/p R)^p\subset
R/pR$ coincides with the set of elements annihilated by all derivations
of $R/pR$ over $\Z/p\Z$ because $R/pR$ is smooth over $\Z/p \Z$. Let us
choose a finite collection
 $$\delta_i\in \Der(R),\,\,i\in
I,\,|I|<\infty$$
 of derivations of $R$ over $\Z$ such that for all
 sufficiently large $p$, the elements $\delta_i\pmod{p}$ span the tangent
 bundle $T_{\spec (R/pR) /\spec (\Z/p\Z)}$.  We have to prove that
 $$\delta_i(f_p^{centr}(y_j))=0\in C_{n,R/pR}\subset
 A_{n,R/pR}\label{form:d=0}$$
 for all $i\in I, j\in\{1,\dots,2n\}$ and almost all prime numbers $p$.
 Let us (for given $i,j,p$) introduce the following notation:
 $$a:=f_p(\xh_j),\,\,b:=\delta_i(f_p(\xh_j))=\delta_i(a)\,\,.$$

 Applying the Leibniz rule to the last expression in the next line
 $$
 \delta_i(f_p^{centr}(y_j))=
 \delta_i(f_p(\xh_j^p))=\delta_i((f_p(\xh_j))^p)=\delta_i(a^p)\,\,,
 $$
we conclude that we have to prove the equality
$$ba^{p-1}+aba^{p-2}+\dots+a^{p-1}b=0\,\,.$$
Notice that $\xh_j$ is a locally ad-nilpotent element of $A_{n,R}$, i.e.
for any $u\in A_{n,R}$ there exists $D=D(u)>0$ such that
$$(\ad (\xh_j))^k(u)=0$$
for $k\ge D(u)$. Namely, one can take $D(u)=\deg(u)+1$ where $\deg(u)$ is the degree of $u$ in the Bernstein
filtration. Using the assumption that $f$ is an automorphism, we conclude 
 that $f(\xh_j)$ is again a locally ad-nilpotent element.
In particular, there exists an integer $D\ge 0$ such that

      $$\left(\ad(f(\xh_j))\right)^D\left(\delta_i(f(\xh_j))\right)=0$$
  for all $i,j$.

Finally, if the prime $p$ is sufficiently large, $p-1\ge D$, then
$$0=(\ad(a))^{p-1}(b)=\sum_{i=0}^{p-1} (-1)^i \binom{p-1}i a^i b \,a^{p-1-i}=\sum_{i=0}^{p-1}a^i b \,a^{p-1-i}\pmod{p}\,\,.$$
This finishes the proof. $\Box$

The conclusion is that for a finitely generated  algebra $R$ smooth over $\Z$ there exists a unique homomorphism
$$\phi_R: \uAut(A_n)(R)\to \uAut(P_n)(R_\infty)$$
such that $\psi_R=\Fr_*\circ \phi_R$. Here $\Fr_*:\uAut(P_n)(R_\infty)\to \uAut(P_n)(R_\infty)$
 is the group homomorphism induced by the endomorphism $\Fr:R_\infty\to R_\infty$ of the coefficient ring.

\begin{conj} \label{conj:lift} In the above notation the image of $\phi_R$ belongs to
$$\uAut(P_n)(i(R)\otimes \Q)\,\,,$$
where $i:R\to R_\infty$ is the tautological inclusion (see Section \ref{sec:rinfty}).
In other words, there exists a unique homomorphism
$$\phi_R^{can}:\uAut(P_n)(R)\to \uAut(P_n)(R\otimes \Q)$$
 such that $\psi_R=\Fr_*\circ i_*\circ \phi_R^{can}$.
\end{conj}

If we assume that the above conjecture holds then we can define a
          constructible map $$\phi^{can}_{n,N}:\uAut^{\le
N}(A_{n,\Q})\to \uAut^{\le N}(P_{n,\Q})$$ for arbitrary integers
$n,N\ge 1$ in the following way.  Let us decompose the scheme $\uAut^{\le
 N}(A_{n,\Q})$ into a finite union $\sqcup_{i\in I} \spec
 R_i,\,\,|I|<\infty$ of closed affine schemes of finite type smooth
over $\spec Q$. Then choose a model $\spec R'_i$ smooth over $\Z$ of
each scheme $\spec R_i$ such that $R'_i$ is a finitely generated ring.
The universal automorphism $u_i$ of $A_n$ defined over $\spec R'_i$
maps under $\phi_R^{can}$ to a certain element $v_i$ of $\uAut^{\le N}
(P_n)(R_i)$.  Taking the union of $v_i,\,i\in I$ we obtain
a constructible map $\phi^{can}_{n,N}$. It is easy to see that this map
does not depend on  the choices made, and the correspondence
$(n,N)\mapsto \phi^{can}_{n,N}$ is compatible with inclusions in
the indices $n,N$. Moreover, the limiting map $\phi^{can}_n$ is compatible
with the group structure.  This is the (conjectured) {\it canonical}
isomorphism between the two automorphism groups.

    \subsection{Continuous constructible maps}

If we do not assume Conjecture \ref{conj:lift}, still the results
of the previous section imply that for any $n,N\ge 1$ there exists
$p_0(n,N)\ge 1$ such that for any prime $p>p_0(n,N)$ we have a
canonically defined constructible map
$$
\phi_{n,N,p}^{can}:\uAut^{\le N}(A_{n,\Z/p\Z} )\hookrightarrow
\uAut^{\le N}(P_{n,\Z/p\Z} )  $$
defined by the property
$$
\Fr_*\circ \,\phi_{n,N,p}^{can} (f)= f^{centr},\,\,\forall f\in \Aut^{\le
N}(A_n)(\k)
$$
for any field $\k$ with $\ch(\k)=p$.

This map is an embedding  because of the following lemma:

\begin{lmm} \label{lmm:incl}
For any field $\k$ of characteristic $p$ the map
$$\Aut(A_{n,\k})\to \Aut(\A^{2n}_\k),\,\,\,f\mapsto f^{centr}$$
 is an inclusion.
\end{lmm}

{\bf Proof}: The above map is a group homomorphism, hence it is enough to
prove that any element $f\in \Aut(A_{n,\k})$ which is mapped to the
identity map is the identity itself.  Let us assume that
$f^{centr}=\Id_{A^{2n}_\k}$ and $f\ne \Id_{A_{n,\k}}$.  We consider two
 cases.

 {\bf First case}. There exists $N\ge 2$ such that $f(\xh_i)$ has
 degree $\le N$ for all $i$, and equal to $N$ for some $i$.  In this
case the same will hold for $f^{centr}$, as the correspondence
$f\mapsto f^{centr}$ preserves the filtration by degree (see Lemma
\ref{lmm:deg}) and is equal to the Frobenius map on the principal
symbols with respect to the  filtration. Hence we get a contradiction
with the assumption $f^{centr}=\Id_{A^{2n}_\k}$.

{\bf Second case}. The degree of $f(\xh_i)$ is equal to $1$ for all
$i,\,\,\,1\le i\le 2n$.  In this case $f$ is an affine symplectic map,
and a direct calculation shows that the corresponding map $f^{centr}$
is also affine symplectic with coefficients equal to the $p$-th power of
those of $f$ (it follows immediately from the results of the next
section).  Hence $f^{centr}=\Id_{A^{2n}_\k}$ implies $f= \Id_{A_{n,\k}}$,
and we again get a contradiction. $\Box$

\begin{conj} \label{conj:isop}
For any $n,N$ there exists $p_1(n,N)\ge p_0(n,N)$ such that
for any prime $p>p_1(n,N)$, the constructible map
$\phi_{n,N,p}^{can}$ is a bijection.
\end{conj}

Obviously, Conjectures \ref{conj:lift} and \ref{conj:isop} together
imply Conjecture \ref{conj:const}.

It is easy to see that the map $\phi_{n,N,p}^{can}$ is {\it continuous} for the
Zariski topology. It follows immediately from the fact that the
correspondence $f\mapsto f^{centr}$ is a regular map (hence continuous)
and that the Frobenius endomorphism of any scheme of finite type in
characteristic $p>0$ is a homeomorphism.  It leads to a natural
question whether the hypothetical canonical isomorphism
$\phi_{n,N}^{can}$ is in fact a homeomorphism for the Zariski topology.

There exists a general notion of  seminormalization ${}^+ S$ for a
reduced scheme $S$ (see \cite{Swan}). One of possible definitions
(in the affine Noetherian case) is that a function $f$ on  ${}^+ S$ is a
reduced closed subscheme $Z_f$ of $S\times \A^1$ which projects
bijectively to $S$. Seminormalization is a tautological operation for
smooth $S$, it coincides with the normalization for integral $S$. The
above question about the homeomorphicity  
 of $\phi_{n,N,p}^{can}$ (a strengthening of Conjecture \ref{conj:isop}) can be rephrased as follows:

{\it   Are  seminormalizations of the {\em reduced}
schemes $\left(\uAut^{\le N}(A_{n,\Z/p\Z} )\right)^{red}$ and
$\left(\uAut^{\le N}(P_{n,\Z/p\Z} )\right)^{red}$  isomorphic? }

\section{Correspondence for tame automorphisms} \label{sec:tame2}

Here we give a proof of Theorem \ref{thm:tame}. Moreover, we will show
that the map $\phi^{can}_n$ is well defined on the tame automorphisms of
 $A_{n,\C}$,  and it takes values in the group of tame automorphisms of
 $P_{n,\C}$.

First of all, we calculate the action of elementary tame automorphisms
 of the Weyl algebra in finite characteristic on its center.

\begin{prop} \label{prop:lin}
Let $\k$ be a field of characteristic $p>2$ and $f\in \Aut(A_{n,\k})$
 be an automorphism given by a linear symplectic mapping on generators

$$f(\xh_i)=\sum_{j=1}^{2n} a_{ij} \xh_j,\,\,\,a_{ij}\in \k\,\,.$$
Then the corresponding automorphism of the center $C_{n,\k}\simeq
 \k[y_1,\dots,y_{2n}]=\k[\xh_1^p,\dots,\xh_{2n}^p]$ is given by
 $$f^{centr}(y_i)=\sum_{j=1}^{2n} (a_{ij})^p y_j\,\,.$$
 \end{prop}

{\bf Proof}: The symplectic group $\Sp(2n,\k)$ is generated by transvections
$$\xh_1\mapsto \xh_1+a\xh_{n+1},\,\,a\in \k,\,\,\,\,\xi_i\mapsto \xi_i
\mbox{ for } i\ge 2$$ and by the Weyl group (Coxeter group $C_n$). The
 correspondence is obvious for elements of the Weyl group, and follows
  from the next Proposition for generalized transvections applied to a
  polynomial of degree $2$.  $\Box$

 \begin{prop} \label{prop:trans} Let $\k$ be a field of characteristic
 $p$ and $f=T_F^A\in \Aut(A_{n,\k})$ be an automorphism corresponding
  to the polynomial $F\in \k[x_1,\dots,x_n]$ (as in Section
 \ref{sec:tame}):
 $$
  T_F^A(\xh_i)=\xh_i,\,\,T_F^A(\xh_{n+i})=\xh_{n+i}+\partial_i F
 (\xh_1,\dots,\xh_n),\,\,1\le i\le n\,\,.
  $$
  Then  one has
 $$f^{centr}(y_i)=y_i,\,\,f^{centr}(y_{n+i})=y_{n+i}+\Fr_*(\partial_i
 F)(y_1,\dots,y_n),\,\,1\le i\le n\,\,,$$
  where the polynomial
 $\Fr_*(\partial_i F)\in \k[y_1,\dots,y_n]$ is obtained from
 $\partial_i F$ by raising all coefficients to the $p$-th power and
  by replacing the variable $x_j$ by $y_j,\,\,1\le j\le n$.  \end{prop}

  {\bf Proof}:
  Let us prove the following identity for the case of one variable:
  $$\left(\frac{d}{dx} +g'\right)^p=\left(\frac{d}{dx}\right)^p
  +\left(g'\right)^p \pmod{p}\,\,,$$
    where $R$ is an arbitrary ring over $\Z/p\Z$, and $g\in R[x]$ is
  any polynomial, $g':={dg}/{dx}$.

  A straightforward calculation over $\Z$ (replace the prime $p$ by an
  integer and use induction) gives

   $$\left(\frac{d}{dx} +g'\right)^p=\sum_{\stackrel{i\ge
  0;\,a_1,\dots,a_{p}\ge 0}{i+\sum_j j a_j=p }}
   \,\frac{p!}{i!\prod_{j=1}^{p} \left(j!\right)^{a_j}\prod_{j=1}^{p}
   a_j!} \,\prod_{j=1}^{p}
   \left(g^{(j)}\right)^{a_j}\left(\frac{d}{dx}\right)^i\,\,.$$

 Here $g^{(j)}$ denotes the $j$-th derivative of the polynomial $g$.
All the coefficients above are divisible by $p$, except for $3$ terms:

$$\left(\frac{d}{dx} +g'\right)^p= \left(\frac{d}{dx}\right)^p +
   (g')^p+  g^{(p)}\pmod{p}\,\,.$$
The last term vanishes because $g^{(p)}=(d/dx)^p(g)=0$ in
characteristic $p$.

For a given $i,\,\,1\le i\le n$ we apply the above identity to
$$R:=\k[\xh_1,\dots,\xh_{i-1},\xh_{i+1},\dots,\xh_n],\,\,\,
g(x):=F(\xh_1,\dots,\xh_{i-1},x,\xh_{i+1},\dots,\xh_n)\,,$$
and get the statement of the Proposition. $\Box$

It follows immediately from Propositions \ref{prop:lin} and
\ref{prop:trans} that $\phi^{can}$ is well defined on symplectic
linear transformations and  transvections $T_F^A$, and hence on all
tame automorphisms of $A_{n,\C}$.  Also it is clear that the following
inclusion holds
$$\Ker \rho^A_n\subset \Ker \rho^P_n$$
(with the notation introduced in Section \ref{sec:tame}). Namely, let us
assume that the composition of a sequence of elementary tame
automorphisms of $A_{n,\C}$ is the identity morphism. The corresponding
transformation of centers of the Weyl algebra in large finite
characteristics is  the  composition (twisted by Frobenius) of the {\it
same} elementary tame automorphisms applied to $P_n$. Hence, the
composition of elementary tame automorphisms of $P_{n,\C}$ is an
identity.

Conversely, let us assume that a composition of elementary
tame transformations of $P_{n,\C}$ is not an identity.  Then applying
Lemma \ref{lmm:incl} we obtain that the corresponding
composition in $\Aut(A_{n,\C})$ is not an identity.  Thus,
$$\Ker \rho^P_n\subset \Ker \rho^A_n\,\,.$$
Theorem \ref{thm:tame} is proven. $\Box$

It is an interesting challenge to find a different proof of Theorem
\ref{thm:tame}, without arguments coming from finite characteristic.

	\section{Conjecture for the inverse map}

Up to now we  talked only  about a homomorphism from $\Aut(A_{n,\C})$
to $\Aut(P_{n,\C})$, and never about the inverse map. Here we
propose a hypothetical construction which produce an automorphism of
the Weyl algebra starting from  a polynomial symplectomorphism.

	 \subsection{Brauer group and $1$-forms}

It is well known that for any Noetherian scheme $S$ in characteristic
$p>0$ there exists a canonical map
$$\alpha:\Omega^1_{abs}(S)/d
\OO(S)\to \Br(S)$$
where $\Omega^1_{abs}(S):=\Gamma(S,\Omega_{S/\spec
\Z/p \Z})$ is the space of global absolute K\"ahler differentials on
$S$. Let us assume for simplicity that $S$ is affine.  For any two
functions $f,g\in \OO(S)$ we define an associative algebra $A_{f,g}$
over $S$ by generators and relations

 $$A_{f,g}:=\OO(S)\langle \xi, \eta\rangle/\left(\mbox{ relations
}[\xi,\eta]=1,\,\, \xi^p=f,\,\,\eta^p=g\right)\,\,.$$

It is easy to see that $A_{f,g}$ is an Azumaya algebra of rank $p$.
The correspondence $\alpha$ is given by
$$\alpha\left(\sum_i f_i d g_i\right):=\sum_i
\left[A_{f_i,g_i}\right]=\left[\otimes_i A_{f_i,g_i}\right]\,\,,$$
where $\left[A_{f_i,g_i}\right]\in \Br(S)$ is the class of algebra
$A_{f_i,g_i}$ in the Brauer group, which is by definition the set of
equivalence classes of Azumaya algebras over $S$
modulo Morita equivalences identical over
centers $\simeq\OO(S)$.

It follows directly from the definitions that for any commutative ring
$R$ over $\Z/ p\Z$ one has an isomorphism of algebras over
 $C_{n,R}=\Center(A_{n,R})\simeq R[y_1,\dots,y_{2n}]$:
$$A_{n,R}\simeq \otimes_{C_{n,R}} A_{y_i,y_{n+i}}\,\,.$$
Hence the class $[A_{n,R}]\in \Br(C_{n,R})$ is given by $1$-form
$$\beta_n:=\sum_{i=1}^n y_i d y_{n+i} \pmod {d C_{n,R}}\,\,.$$

 The correctness of the definition of the map $\alpha$ follows
from the existence of certain bimodules establishing the Morita
equivalences.

\subsubsection{ Explicit Morita equivalences}

One can construct explicitly the following  isomorphisms of $\OO(S)$-algebras:
\begin{itemize} 
\item $A_{f,0}\simeq A_{0,g}\simeq \Mat(p\times p,\OO(S))$
\item  $A_{f,g}\simeq A_{g,-f}$  (Fourier transform)
\item $A_{f_1+f_2,g}\otimes_{\OO(S)}  \Mat(p\times p,
\OO(S))\simeq  A_{f_1,g}\otimes_{\OO(S)} A_{f_2,g}$
\item $A_{f,g_1+g_2}\otimes_{\OO(S)}  \Mat(p\times p,\OO(S))\simeq  
A_{f,g_1}\otimes_{\OO(S)} A_{f,g_2}$
\item $A_{f,gh}\otimes_{\OO(S)} A_{g,hf}\otimes_{\OO(S)}
 A_{fh,g}\simeq Mat(p^{3}\times p^{3},\OO(S))$
\item $A_{1,f}\simeq \Mat(p\times p,\OO(S))$
\end{itemize}
  corresponding to basic identities in $\Omega^1_{abs}(S)/d \OO(S)$:
 \begin{itemize}
 \item $ f \,d(0)= 0\,dg =0$
 \item $f \,dg=g\, d(-f) \in \Omega^1_{abs}(S)/d \OO(S)$
 \item $ (f_1+f_2)\,dg=f_1 dg +f_2 dg$
 \item $  f \,d(g_1+g_2)=f\,dg_1+f\,dg_2$
 \item $  f\,d(gh)+g\,d(hf)+ h\,d(fg)=0$
 \item $ 1\,df=0\in \Omega^1_{abs}(S)/d \OO(S)$
 \end{itemize} 
 
 It is convenient to replace the matrix algebra $\Mat(p\times p,\OO(S))$ by the algebra $A_{0,0}$.
 For example, the isomorphism 
  between the algebras $A_{0,0}$ and $A_{f,0}$,
 $\xi^p=0,\eta^p=0, [\xi,\eta]=1$ and
 $\xi'^p=f,\eta'{}^p=0, [\xi',\eta']=1$ corresponding to the pairs $(f,0)$ and
 $(0,0)$, is given by the formula $\xi'\to\xi-f\eta^{p-1}$, $\eta'\to\eta$.
 Similarly,  the more complicated fifth isomorphism in the above list   
 $$
 A_{f,gh}\otimes_{\OO(S)}
A_{g,hf}\otimes_{\OO(S)}A_{h,fg}\simeq
 A_{f,0}\otimes_{\OO(S)}
A_{g,0}\otimes_{\OO(S)}A_{h,0}(\simeq
\Mat(p^{3n}\times
p^{3n},\OO(S)))\,\,
 $$
(here we use isomorphisms $A_{0,0}\simeq A_{f,0}$ etc.) is given by the formula
 $$
 \xi_1'=\xi_1,\ \eta_1'=\eta_1-\xi_2\xi_3\,,\
 \xi_2'=\xi_2,\ \eta_2'=\eta_2-\xi_3\xi_1\,,\
 \xi_3'=\xi_3,\ \eta_3'=\eta_3-\xi_1\xi_2\,.
 $$
We leave to the interested reader the construction of other isomorphisms as an exercise.

    \subsection{Pullback of the Azumaya algebra under a symplectomorphism}

    Let $R$ be a finitely generated smooth commutative algebra over
$\Z$,  and $g\in \Aut(P_{n,R})$ be a symplectomorphism defined over
$R$. Our goal is to construct an automorphism $f\in \Aut(A_{n,R\otimes
Q})$ such that $\phi_R^{can}(f)=g$.

By definition, we have $$g(\omega)=\omega$$ where $\omega
=d\beta_n=\sum_{i=1}^n dx_i\wedge dx_{i+n}$ is the standard symplectic
form on $\A^{2n}_R$.  Hence there exists a polynomial $P\in \Q\otimes
 R[x_1,\dots,x_{2n}]$ such that

  $$g(\beta_n)=\beta_n+d P \in \Omega^1_{\Q\otimes R[x_1,\dots,x_{2n}]/\Q\otimes R}\,,$$
the reason is that $H^1_{de\,\,Rham}(\A^{2n}_{R\otimes \Q})=0$.  We
  can add to $R$ inverses of finitely many primes and assume that $P\in
  R[x_1,\dots,x_{2n}]$.

For any prime $p$ let us consider the symplectomorphism $\Fr_*(g)\in
\Aut(P_{n,R/pR})$ which is obtained from $g_p:=g\pmod{p}\in
\Aut(P_{n,R/pR})$ by raising to the $p$-th power all the coefficients (in
other words, by applying the Frobenius endomorphism $R/pR\to
R/pR$).

We claim that $\Fr_*(g_p)$ preserves the class of $\beta_n$ in
$$\Omega^1_{abs}\left(R/p R[x_1,\dots,x_{2n}]\right)/ d \left(R/p
R[x_1,\dots,x_{2n}]\right)\,\,.$$

Obviously, we have an identity

$$\Fr_*(g_p)(\beta_{n,p})=\beta_{n,p}+ d\,\Fr_*(P_{p})\,\in
\Omega^1_{R/pR[x_1,\dots,x_{2n}]/(R/pR)}$$
 in the space of {\it relative} 
 1-forms over $R/pR$, where
 $\beta_{n,p}:=\beta_n\pmod{p},\,\,P_p:=P\pmod{p}$.  The same identity
 holds in  {\it absolute} 1-forms because all the coefficients of
  the transformation $\Fr_*(g_p)$ and of the polynomial $\Fr_*(P_p)$ belong to
the image of the Frobenius map $\Fr(R/pR)=(R/pR)^p\subset R/pR$, and hence
behave like constants for absolute $1$-forms:
 $$d (a^p b)= a^p db \in
 \Omega^1_{abs} \left(R/p R[x_1,\dots,x_{2n}]\right), \,\,\,\forall
 a\in R/pR,\,\,b\in \left(R/p R[x_1,\dots,x_{2n}]\right)\,\,.$$

The conclusion is that the pullback of the algebra $A_{n,R/pR}$ under
the symplectomorphism $\Fr_*(g_p)$ has the same class in
$\Br(\A^{2n}_{R/pR})$ as the algebra $A_{n,R/pR}$ itself.  Therefore, there
exists a Morita equivalence between these two algebras, identical on the
center.  In other words, we proved that there exists a Morita
autoequivalence of $A_{n,R/pR}$ inducing an automorphism $\Fr_*(g_p)$ of
the center.

Let us denote by $M_{g,p}$ {\it any} bimodule over $A_{n,R/pR}$
corresponding to the above Morita autoequivalence. The following result
shows that this bimodule is essentially unique, its isomorphism class
is uniquely determined by $g$.

\begin{lmm}\label{lmm:morita}
Any bimodule over $A_{n,R/pR}$ inducing Morita autoequivalence
identical on the center $C_{n,R/pR}$ is isomorphic to the diagonal
bimodule.  Any automorphism of the diagonal bimodule is given by
multiplication by a constant in $(R/pR)^\times $.
\end{lmm}

 {\bf Proof}: It is easy to see that isomorphism classes of such
bimodules form a torsor over $H^1_{\acute{e} t}(\A^{2n}_{R/pR},
\Gm)\simeq 0$.  Similarly, symmetries of any such modules are
   $H^0_{\acute{e} t}(\A^{2n}_{R/pR}, \Gm)\simeq (R/pR)^\times$. $\Box$

    \subsection{A reformulation of the conjectures}

In the notation introduced above, the bimodule $M_{g,p}$ is a finitely
generated projective left $A_{n,R/pR}$-module.  It corresponds to an
automorphism of the algebra $A_{n,R/pR}$ iff it is a free rank one module.
Moreover, the corresponding automorphism is uniquely defined because
all  invertible elements of $A_{n,R/pR}$ are central.

Conjecture \ref{conj:isop} is equivalent to the following

\begin{conj}\label{conj:isop2} For any finitely generated smooth commutative algebra $R$ over $\Z$ and any
$g\in \Aut(P_{n,R})$ for all sufficiently large $p$, the bimodule $M_{g,p}$ is a free rank one left $A_{n,R/pR}$-module.
\end{conj}

There is no clear evidence for this conjecture as there are examples
 (see \cite{Stafford}) of projective finitely generated modules over
the Weyl algebra $A_1$ in characteristic zero, which are in a sense of
rank $1$ and not free.  In other words, an analogue of the Serre conjecture
for Weyl algebras is false.

 If Conjectures \ref{conj:isop} and \ref{conj:isop2} fail to be true, it is still quite feasible that
  the following {\it weaker} version of Conjecture \ref{conj:iso} holds:

\begin{conj} \label{conj:morita}
The group of Morita autoequivalences of the algebra $A_{n,\C}$ is
isomorphic to the group of polynomial symplectomorphisms
 $Aut(P_{n,\C})$.
 \end{conj}

Its equivalence to Conjecture \ref{conj:iso} depends on the answer to
the following question:

{\it Does any Morita autoequivalence of  $A_{n,\C}$ come from an
automorphism?}

     \subsection{Adding the Planck constant}

It follows from  Conjecture \ref{conj:iso}  that there should exist
 a mysterious nontrivial action of $\C^\times$
by {\it outer} automorphisms of the group $\Aut(A_{n,\C})$:
    $$\C^\times\to \mbox{Out}(\Aut(A_{n,\C}))$$
corresponding to the conjugation by the dilations
$$(x_1,\dots,x_{2n})\to (x_1,\dots,x_n,\lambda x_{n+1},\dots,\lambda
x_{2n}),\,\,\lambda\in \C^\times$$
acting by automorphisms of $\Aut(P_{n,\C})\subset \Aut(\C[x_1,\dots,
x_n])=\Aut(\A^{2n}_\C)$.  Alternatively, one can say that there is a
one-parameter family of hypothetical isomorphisms between
$\Aut(A_{n,\C})$ and $\Aut(P_{n,\C})$.

In general, it makes sense to introduce a new central variable $\hbar$
  (``Planck constant'' parameterizing the above family of
isomorphisms), and define the algebra $A_{n,R}^{\hbar}$ as an associative
algebra over the commutative ring $R[\hbar]$ given by generators
$\xh_1,\dots,\xh_{2n}$ and defining relations
$$[\xh_i,\xh_j]=\hbar \omega_{ij}\,\,.$$
We propose the following

\begin{conj}\label{conj:final}
For any finitely generated smooth commutative algebra $R$ over $\Z$
and any symplectomorphism $g\in \Aut(P_{n,R})$ there exists a positive
integer $M$ and an automorphism $\tilde{g} \in
\Aut(A^{\hbar}_{n,R(M^{-1})[\hbar]})$ over $R(M^{-1})[\hbar]$ such that
 \begin{itemize}
 \item $\tilde{g} \pmod \hbar = g$,
  \item for all sufficiently large $p$ the automorphism $\tilde{g}\pmod p$
preserves the subalgebra $R/pR[ y_1^p,\dots, y_{2n}^p]$.
  \end{itemize}
  \end{conj}

This conjecture seems to be the best one, one can  easily see that it
  implies {\it all} the conjectures previously made in the present paper.

     \section{On the extensions of the conjecture to other algebras}

The Weyl algebra is isomorphic to the algebra of differential operators on
the affine space.  It has a natural generalization, the algebra $\D(X)$
of differential operators on a smooth affine algebraic variety
$X/\k,\,\,\ch(\k)=0$. The corresponding Poisson counterpart is the
algebra of functions on the cotangent bundle $T^* X$ endowed with the
natural Poisson bracket.  One may ask what happens with our conjectures for
such algebras. Unfortunately, one of our key results (Theorem
\ref{thm:frob}) relies heavily on the local ad-nilpotence property of the generators of 
$A_{n,\k}=\D(\A^n_\k)$ which does not hold in general.  Moreover, it is
easy to see that there are counter-examples to the naive extension of
 Conjecture \ref{conj:iso}.  In particular, for
$X=\A^1_\k\setminus\{0\}$ with invertible coordinate $x$, the
automorphism of $D(X)$ given by

 $$x\mapsto x,\,\,\,\partial/\partial x\mapsto \partial/\partial x
 +c/x$$
does not seem to correspond to any particular symplectomorphism of
 $T^*X$, as the corresponding transformation of the center in
characteristic $p>0$ is

  $$x^p\mapsto x^p,\,\,(\partial/\partial x)^p\mapsto
(\partial/\partial x)^p +(c^p-c)/x^p\,\,.$$
  The constant $(c^p-c)$ does not belong to the image of the Frobenius map
  in general.

One can try to generalize Conjecture \ref{conj:iso} in a different
  direction.  Any automorphism $f$ of $A_{n,\k}$ gives a bimodule,
which can be interpreted as a holonomic module $M_{(f)}$ over
$A_{2n,\k}$. Similarly, any symplectomorphism $g\in \Aut(P_{n,\k})$
gives a Lagrangian submanifold $L_{(g)}\subset \A^{4n}_\k$ (the graph
 of $g$). The idea is to establish a correspondence between holonomic
 $\D$-modules and Lagrangian subvarieties (a version of the Hitchin
 correspondence).

In general, any holonomic module $M$ over $D(X)$ defined over a
finitely generated smooth ring $R\subset \C$, gives a family of
coherent sheaves $M_p$ over $\Fr_p^*(T^*X)$ by considering reductions
modulo $p$ as modules over the centers of algebras of differential
operators. One expects (in analogy with the theory of characteristic
varieties) that the support of $M_p$ is a Lagrangian subvariety $L_p$
of the twisted by Frobenius cotangent space, at least for large $p$
(see the recent preprint \cite{OV} where the analog of the Hitchin
correspondence was studied for a {\it fixed} prime $p$).  The subvariety
$L_p$ can be singular and not necessarily conical.  Moreover, its
dependence on $p$ in general seems to be chaotic.

Nevertheless, we expect that in certain circumstances there is a
  canonical correspondence in characteristic zero.  Namely, we believe
that for any closed Lagrangian subvariety  $L$ of $T^*X$ such that $L$
 is smooth and $H^1(L(\C),\Z)=0$, there exists a canonical
holonomic $D(X)$-module $M_L$ such that $(M_L)_p$ is supported on
$\Fr_p(L)$ and moreover, is locally isomorphic to the sum of
$p^{dim\,X}$ copies of $\OO_{\Fr_p(L)}$. This would imply Conjecture
\ref{conj:morita}.

\vspace{5mm}

Addresses:

A.B.-K.: Institute of Mathematics, Hebrew University, Givat Ram, Jerusalem
91904,
 Israel.

{kanel@mccme.ru} \\

M.K.: IHES, 35 route de Chartres, Bures-sur-Yvette 91440, France

{maxim@ihes.fr}

\end{document}